\documentclass{amsart}

\usepackage[latin1]{inputenc}
\usepackage[english]{babel}
\usepackage{amsmath}
\usepackage{amssymb}
\usepackage{centernot}
\usepackage{caption}
\usepackage{color}
\usepackage{subcaption}
\usepackage{etoolbox}
\usepackage{upgreek}
\usepackage{mathabx}
\usepackage{titlesec}
\usepackage{enumitem}

\patchcmd{\section}{\scshape}{\bfseries}{}{}
\makeatletter
\renewcommand{\@secnumfont}{\bfseries}
\makeatother
\newcommand*{\justifyheading}{\raggedright}
\titleformat {\section}
  {\normalfont \Large \bfseries \justifyheading}{\thesection.}{1em}{}
\titleformat {\subsection}
  {\normalfont \large \bfseries \justifyheading}{\thesubsection.}{1em}{}
\patchcmd{\abstract}{\scshape\abstractname}{\textbf{\abstractname}}{}{}

\numberwithin{equation}{section}
\newtheorem{thm}{Theorem}[section]
\newtheorem{prop}[thm]{Proposition}
\newtheorem{lem}[thm]{Lemma}
\newtheorem{cor}[thm]{Corollary}

\newtheorem{obs}[thm]{Observation}
\newtheorem{theo}[thm]{Theorem}
\newtheorem{defn}[thm]{Definition}
\newtheorem{nota}[thm]{Notation}

\theoremstyle{definition}

\newtheorem{rem}[thm]{Remark}

\newenvironment{proof2}[1][Proof.]{\begin{trivlist}
\item[\hskip \labelsep {\bfseries \itshape #1}]}{\end{trivlist}}
\newenvironment{proof3}[1][Proof.]{\begin{trivlist}
\item[\hskip \labelsep {\itshape #1}]}{\end{trivlist}}

\newcommand{\Arr}[1]{\mathcal{A}_{#1}}

\newcommand{\R}{\mathbb{R}}
\newcommand{\N}{\mathbb{R}_{\geq 0}}
\newcommand{\PP}{\mathbb{P}}

\newcommand{\card}[1]{\left|#1\right|}

\newcommand{\Frac}[2]{\displaystyle\frac{#1}{#2}}

\newcommand{\eigspace}[1]{\mathbf{E}_{\uplambda_{#1}}}

\usepackage{graphicx}
\usepackage{mathrsfs}
\usepackage[foot]{amsaddr}

\usepackage[round]{natbib}

\title[Largest eigenvalue of the Laplacian]{Largest eigenvalue of the Laplacian matrix}
\author{Benjamin Iriarte}
\address{Department of Mathematics, Massachusetts Institute of Technology, Cambridge MA, 02139, USA}
\email{biriarte@math.mit.edu}


\begin{document}
\begin{abstract}
We study the eigenspace of the Laplacian matrix of a simple graph corresponding to the largest eigenvalue, subsequently arriving at 
the theory of modular decomposition of T. Gallai.      
\end{abstract}

\maketitle

\section{Introduction.}\label{sec:intro}

Let $G=G([n],E)$ be a simple (undirected) graph, where $[n]=\{1,2,\dots,n\}$, $n\in\PP$. The \emph{adjacency matrix of $G$} is the 
$n\times n$ matrix $A=A(G)$ such that:
$$(A)_{ij}=a_{ij}:=\left\{
\begin{array}{ll}
1&\text{if $\{i,j\}\in E$,}\\
0&\text{otherwise.}\end{array}\right.$$ 
The \emph{Laplacian matrix of $G$} is the $n\times n$ matrix $L=L(G)$ such that:
$$(L)_{ij}=l_{ij}:=\left\{
\begin{array}{ll}
d_i&\text{if $i=j$,}\\
-a_{ij}&\text{otherwise,}\end{array}\right.$$ 
where $(d_{G})_i=d_i$ is the degree of vertex $i$ in $G$. 

The spectral theory of these matrices, {\it i.e.} the theory about their eigenvalues and eigenspaces, has been the object of much study for the last $40$ years. The roots of this
beautiful theory, however, can arguably be traced back to \emph{Kirchhoff's matrix-tree theorem}, whose first proof is often attributed to~\cite{borchardt} even though at least one
proof was already known by~\cite{sylvester}. A recollection of some interesting applications of the theory can be found in~\cite{spiel}, and more complete accounts of the mathematical backbone are~\cite{brouwer} and~\cite{chung}. Still, it would be largely inconvenient and prone to unfair omissions to attempt here a fair account
of the many contributors and contributing papers that helped shape the state-of-the-art 
of our knowledge of graph spectra, and we refer the reader to our references for further inquiries of the literature. 

This article aims to fill one (of the many) gap (s) in our current knowledge of the theory, namely, the lack of results about eigenvectors of the Laplacian with largest eigenvalue. We will answer the question: What information about the structure of a graph is carried in these eigenvectors?
Our work follows the spirit of~\cite{fiedler1}, who pioneered the use of eigenvectors of the Laplacian matrix to learn about a graph's structure. One of the first observations that can be made about $L$ is that it is positive-semidefinite, a consequence of it being a product of \emph{incidence matrices}. We will thus let, 
$$0=\uplambda_1\leq \uplambda_2\leq \dots\leq\uplambda_n=\uplambda_{\max}=\uplambda_{\max}(G),$$ 
be the (real) eigenvalues of $L$, and note that
$\uplambda_2>0$ if $G$ is a connected graph; we have effectively dropped $G$ from the notation for convenience but remark that eigenvalues and eigenvectors 
depend on the particular graph at question, which will be clear from the context. We will also let $\eigspace{i}$ be the eigenspace corresponding to
$\uplambda_i$. In its most primitive form, \emph{Fiedler's nodal domain's theorem} [\cite{fiedler2}] states that when $G$ is connected and for all $x\in\eigspace{2}$,
the induced subgraph $G\left[\{i\in[n]:x_i\geq 0\}\right]$ is connected. Related work, also relevant to the present writing, might be found in~\cite{merris}.  

We will go even further in the way in which we use eigenvectors of the Laplacian to learn properties of $G$. To explain this, let us firstly call a map, 
$$O:E\rightarrow \left([n]\times [n]\right)\cup E=[n]^2\cup E,$$ 
such that $O(e)\in\{e,(i,j),(j,i)\}$ for all $e:=\{i,j\}\in E$, an \emph{(partial) orientation} of $E$ (or $G$), and say that, furthermore, $O$ is \emph{acyclic} if $O(e)\neq e$ for all $e$ and the directed-graph
on vertex-set $[n]$ and edge-set $O(E)$ has no directed-cycles. On numerous occasions, we will somewhat abusively also identify
$O$ with the set $O(E)$.
 
During this paper, eigenvectors of the Laplacian and more precisely, elements of
$\eigspace{\max}$, will be used to obtain orientations of certain 
(not necessarily induced) subgraphs of $G$. 
Henceforth, given $G$ and
for all $x\in\R^{[n]}$, 
the reader should always 
automatically consider the orientation (map) $O_{x}=O_{x}(G)$ associated to $x$, 
$O_x:E\rightarrow [n]^2\cup E$, such that for $e:=\{i,j\}\in E$: 
\begin{align*}
O_x(e)=\left\{
\begin{array}{ll}
e&\text{if $x_i=x_j$,}\\
(i,j)&\text{if $x_i<x_j$,}\\
(j,i)&\text{if $x_i>x_j$.}
\end{array}
\right.
\end{align*}
The orientation $O_x$ will be said to be \emph{induced} by $x$ ({\it e.g.} Figure~\ref{fig:ex3}). 

Implicit above is another subtle perspective that we will adopt, explicitly, that vectors
$x\in\R^{[n]}$ are real functions from the vertex-set of the graph in question (all our graphs will be on vertex-set $[n]$). In our case, this graph is $G$, 
and even though accustomed to do so otherwise, 
entries of $x$ should be really thought of as
being indexed by vertices of $G$ and not simply by positive integers. Later on in Section~\ref{sec:lecg}, for example, we will regularly
state (combinatorial) results about the \emph{fibers} of $x$ when $x$ belongs to a certain subset of $\R^{[n]}$ ({\it e.g.}
$\eigspace{\max}$), thereby regarding these fibers as vertex-subsets of the particular graph being discussed at that moment. 

Using this perspective, we will learn that the eigenspace $\eigspace{\max}$ is closely related 
to the theory of \emph{modular decomposition} of~\cite{gallai}; orientations induced by elements of $\eigspace{\max}$
lead naturally to the discovery of \emph{modules}. This connection will most concretely be exemplified 
when $G$ is a \emph{comparability graph}, in which case these orientations iteratively
correspond to and exhaust the \emph{transitive orientations} of $G$. It will be instructive to see 
Figure~\ref{fig:ex} at this point.

In Section~\ref{sec:bd}, we will introduce the background and definitions necessary to state   
the precise main contributions of this article. These punch line results will then be presented in Section~\ref{sec:lecg}. The central theme of 
Section~\ref{sec:lecg} will be a stepwise proof of Theorem~\ref{theo:mainlap}, 
our main result for comparability graphs, which summarily states that when $G$ is a comparability graph, elements of $\eigspace{\max}$ induce
transitive orientations of the \emph{copartition subgraph} of $G$. It will be along the natural course of this proof that we 
present our three main results that apply to arbitrary simple graphs: Propositions~\ref{prop:laplamod1} 
and~\ref{prop:laplamod2}, and Corollary~\ref{cor:laplamod}.  

Finally, in Section~\ref{sec:ccg}, we will present a curious novel characterization of comparability graphs that results from the
theory of Section~\ref{sec:lecg}.  
  
\section{Background and definitions.}\label{sec:bd}    
\subsection{The graphical arrangement.}
\begin{defn}\label{defn:ga}
Let $G=G([n],E)$ be a simple (undirected) graph. The \emph{graphical arrangement of $G$} is the union of hyperplanes in $\R^{[n]}$:
$$\Arr{G}:=\{x\in\R^{[n]}:x_i-x_j=0\text{ , $\forall$ $\{i,j\}\in E$}\}.$$
\end{defn}

Basic properties of graphical arrangements and, more generally, of hyperplane arrangements, are presented in Chapter 2 of~\cite{ste3}.

For $G$ as in Definition~\ref{defn:ga}, let $\mathcal{R}(\Arr{G})$ be the collection of all (open) connected components of the set
$\R^{[n]}\backslash \Arr{G}$. An element of $\mathcal{R}(\Arr{G})$ is called a \emph{region} of $\Arr{G}$, and every region of $\Arr{G}$ is therefore an
$n$-dimensional open convex cone in $\R^{[n]}$. Furthermore, the following is
true about regions of the graphical arrangement:
\begin{prop}\label{prop:gaao}
Let $G$ be as in Definition~\ref{defn:ga}. Then, for all $R\in \mathcal{R}(\Arr{G})$ and $x,y\in R$, we have that:
$$O_R:=O_x=O_y.$$ 
Moreover, the map $R\mapsto O_R$ from the set of regions of $\Arr{G}$ to the set of orientations of $E$ 
is a bijection between $\mathcal{R}(\Arr{G})$ and the set of \emph{acyclic orientations} of $G$.
\end{prop}

Motivated by Proposition~\ref{prop:gaao} and the comments before, we will introduce special notation for certain subsets of $\R^{[n]}$ obtained
from $\Arr{G}$. 

\begin{nota}\label{nota:cones}
Let $G$ be as in Definition~\ref{defn:ga}. For an acyclic orientation $O$ of $E$, we will let $C_O$ denote the 
$n$-dimensional closed convex cone in $\R^{[n]}$ that is equal to the topological closure of the region of $\Arr{G}$ 
corresponding to $O$ in Proposition~\ref{prop:gaao}. 
\end{nota}

\subsection{Modular decomposition.}

We need to concur on some standard terminology and notation from graph theory, so 
let $G=G([n],E)$ be a simple (undirected) graph and $X$ a subset of $[n]$. 

As customary, $\overline{G}$ denotes the \emph{complement graph} of $G$. 
The notation $N(X)$ denotes 
the \emph{open neighborhood} of
$X$ in $G$: 
$$N(X):=\left\{j\in[n]\backslash X:\text{ there exists some }i\in X\text{ such that }\{i,j\}\in E\right\}.$$
The \emph{induced
subgraph} of $G$ on $X$ is denoted by $G[X]$, and the binary operation
of graph \emph{disjoint union} is represented by the plus sign $+$. 
Lastly, for $Y\subseteq [n]$, $X$ and $Y$
are said to be \emph{completely adjacent} in $G$ if:
\begin{align*}
\text{$X\cap Y=\emptyset$, and }\\
\text{for all $i\in X$ and $j\in Y$, we have that
$\{i,j\}\in E$.}
\end{align*}

The concepts of \emph{module} and \emph{modular decomposition} in graph theory were introduced by~\cite{gallai} as a means to understand the structure
of comparability graphs. The same work would eventually present a remarkable characterization of these graphs in terms of 
forbidden subgraphs. Section~\ref{sec:lecg} of the present work will present an alternate and surprising route to modules. 

\begin{defn}\label{defn:module}
Let $G=G([n],E)$ be a simple (undirected) graph. A \emph{module} of $G$ is a set $A\subseteq [n]$ such that for all $i,j\in A$: 
$$N(i)\backslash A=N(j)\backslash A=N(A).$$ 
Furthermore, $A$ is said to be \emph{proper} if $A\subsetneq[n]$, \emph{non-trivial} if $\card{A}>1$, and
\emph{connected} if $G[A]$ is connected. 
\end{defn}  
\begin{cor}\label{cor:trivial}
In Definition~\ref{defn:module}, two disjoint modules of $G$ are either completely adjacent or no edges exist between them. 
\end{cor}

Let us now present some basic results about modules that we will need.
\begin{lem}[\cite{gallai}]\label{lem:modbasic}
Let $G=([n],E)$ be a connected graph such that $\overline{G}$ is connected. If $A$ and $B$ are maximal (by inclusion) proper modules
of $G$ with $A\neq B$, then $A\cap B=\emptyset$.
\end{lem}
\begin{cor}[\cite{gallai}]\label{cor:modbasic}
Let $G=([n],E)$ be a connected graph such that $\overline{G}$ is connected. Then, there exists a unique partition of $[n]$ into
maximal proper modules of $G$, and this partition contains more than two blocks. 
\end{cor}

From Corollary~\ref{cor:modbasic}, it is therefore natural to consider the partition
of the vertex-set of a graph into its maximal modules; the appropriate
framework for doing this is presented in Definition~\ref{defn:decomp}. Hereafter, however, we will 
assume that our graphs are connected unless otherwise stated since {\bf (1)} the results for disconnected graphs will
follow immediately from the results for connected graphs, and {\bf (2)} this will 
allow us to focus on the interesting parts of the theory.     

\begin{defn}[\cite{pg}]\label{defn:decomp}
Let $G=G([n],E)$ be a connected graph. We will let
the \emph{canonical partition} of $G$ be the set $\mathcal{P}=\mathcal{P}(G)$ such that:
\begin{itemize}
\item[a.] If $\overline{G}$ is connected, $\mathcal{P}$ is the unique partition of $[n]$ into the maximal proper modules of $G$.
\item[b.] If $\overline{G}$ is disconnected, $\mathcal{P}$ is the partition of $[n]$ into the vertex-sets of the connected components of
$\overline{G}$. 
\end{itemize}
\end{defn}
Hence, in Definition~\ref{defn:decomp}, every element of the
canonical partition is a module of the graph. Elements of the canonical partition of a graph
on vertex-set $[8]$ are shown in Figure~\ref{fig:ex2}.

\begin{defn}\label{defn:cpgraph}
In Definition~\ref{defn:decomp}, we will let the \emph{copartition subgraph} of $G$ be the graph $G^{\mathcal{P}}$ 
on vertex-set $[n]$ and edge-set equal to: 
$$E\big\backslash\left\{\{i,j\}\in E:i,j\in A\text{ for some }A\in\mathcal{P}\right\}.$$ 
\end{defn}
\subsection{Comparability graphs.}

We had anticipated the importance of comparability graphs in this work, yet, we need to define what they are.  

\begin{figure}[ht]
\begin{tabular}{llll}
\begin{subfigure}{.25\textwidth}
  \centering
   \caption{}
  \includegraphics[width=0.9\linewidth]{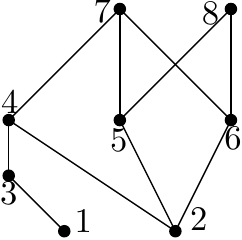}
  \label{fig:ex1}
\end{subfigure}
&
\begin{subfigure}{.25\textwidth}
  \centering
   \caption{}
  \includegraphics[width=0.9\linewidth]{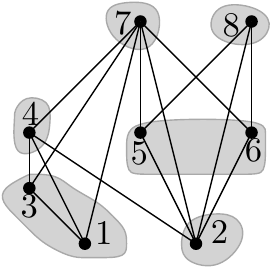}
  \label{fig:ex2}
\end{subfigure}
&
\multicolumn{2}{c}{
\begin{subfigure}{.5\textwidth}
  \centering
   \caption{}
  \includegraphics[width=1\linewidth]{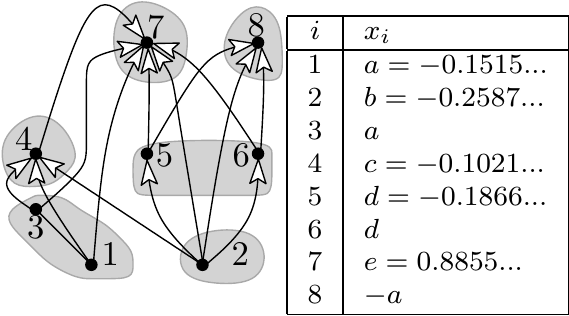}
  \label{fig:ex3}
\end{subfigure}
}
\end{tabular}

\caption{{\bf (A)} Hasse diagram of a poset $P$ on $[8]$. {\bf (B)} Comparability graph $G=G([8],E)$
of the poset $P$, where closed regions are maximal proper modules of $G$. {\bf (C)} Unit
eigenvector $x\in \eigspace{\max}$ of $G$ fully calculated, where $\dim\left\langle \eigspace{\max}\right\rangle = 1$. Arrows represent
the induced orientation $O_x$ of $G$. Notice the relation between $O_x$, modules of $G$, and poset $P$.}
\label{fig:ex}
\end{figure}

\begin{defn}\label{def:compar}
A \emph{comparability graph}
is a simple (undirected) graph $G=G(V,E)$ such that there exists a partial
order on $V$ under which
two different vertices $u,v\in V$ are comparable if and only if
$\{u,v\}\in E$.  
\end{defn}
A comparability graph on vertex-set $[8]$ is shown in Figure~\ref{fig:ex2}.

Comparability graphs are \emph{perfectly orderable graphs} and more
generally, \emph{perfect graphs}. These three families of graphs are all 
large hereditary classes of graphs. 

Note that, given a comparability graph $G=G(V,E)$, we can find at least two
partial
orders on $V$ whose comparability
graphs (obtained as discussed in Definition~\ref{def:compar}) agree precisely with $G$, and the
number of such partial orders depends on the modular decomposition of $G$. Let us
record this idea in a definition.

\begin{defn}\label{def:trans}
Let $G=G(V,E)$ be a comparability graph, and let $O$ be an acyclic
orientation of $E$. Consider the partial order induced by $O$ under which,
for $u,v\in V$, $u$ is less than $v$ iff there is a directed-path in $O$ that begins in
$u$ and ends in $v$. If the comparability graph of this partial order
on $V$ (obtained as in Definition~\ref{def:compar}) agrees precisely with $G$, then we will say that $O$
is a \emph{transitive orientation} of $G$.
\end{defn}

\subsection{Linear algebra.}

Some standard terminology of linear algebra and other related conventions that we adopt are presented here. Firstly, 
we will always be working
in Euclidean space $\R^{[n]}$, and all (Euclidean-normed real) vector spaces considered are assumed to live therein. 
Euclidean norm is denoted by $||\cdot||$.
The standard basis of
$\R^{[n]}$ will be $\{e_i\}_{i\in[n]}$, as customary. Generalizing this notation, for
all $I\subseteq [n]$, we will also let:
$$e_{I}:=\sum_{i\in I}e_i.$$ 
The orthogonal complement
in $\R^{[n]}$ to $\text{span}_\R  \left\langle e_{[n]}\right\rangle$ will be of importance to us, so we will use special notation to denote it:
$$\R^{\ast[n]}:=\left(\text{span}_\R  \left\langle e_{[n]}\right\rangle\right)^{\perp}.$$

For an arbitrary vector space $\mathcal{V}$ and a linear transformation $T:\mathcal{V}\rightarrow\mathcal{V}$, we will say that 
a set $U\subseteq\mathcal{V}$ is \emph{invariant under} $T$, or that $T$ is \emph{$U$-invariant}, if 
$T(U)\subseteq U$. Lastly, a key concept of this paper: 
\begin{align*}
\text{For a vector $x\in\R^{[n]}$ and a set $\upxi \subseteq[n]$,} \\
\text{we will say that $\upxi$ is a \emph{fiber} of $x$}\\
 \text{if
there exists $\upalpha\in\R$ such that $x_i=\upalpha$ if and only if $i\in \upxi$.}
\end{align*}

The notion of being a \emph{generic} vector in a certain vector space, to be understood from the point of view of
\emph{Lebesgue measure theory}, is a central ingredient in many of our results. We now make this notion precise. 

\begin{defn}\label{defn:ucuv}
Let $\mathcal{V}$ be a linear subspace of $\R^{[n]}$ with $\dim\left\langle \mathcal{V}\right\rangle >0$. 
We will say that a vector $x\in \mathcal{V}$ is a \emph{uniformly chosen at random unit vector or 
{\bf u.c.u.v.}} if $x$ is 
uniformly chosen at random from the set $\{y\in \mathcal{V}:||y||=1\}$. 

For $x\in \mathcal{V}$ a \emph{u.c.u.v.}, a certain event or statement about $x$ is said to occur or hold true \emph{almost surely}
if it is true with probability one. 
\end{defn}
\subsection{Spectral theory of the Laplacian.}\label{subsec:stl}

We will need only a few background results on the spectral theory of the Laplacian matrix of a graph. 
We present these below in a single statement, but refer
the reader to \cite{brouwer} for additional background and history. 

\begin{lem}\label{lem:tech1}
Let $G=G([n],E)$ be a simple (undirected) graph. Let $L=L(G)$ be the Laplacian matrix of 
$G$ and $0=\uplambda_1\leq\uplambda_2\leq\dots\leq \uplambda_n=\uplambda_{\max}=\uplambda_{\max}(G)$ be the eigenvalues of
$L$. Then:
\begin{itemize}
\item[{\bf 1.}] The number of connected components of $G$ is equal to the multiplicity of the eigenvalue $0$ in $L$.
\item[{\bf 2.}] If $\overline{G}$ is the complement of $G$ and $\overline{L}$ is the Laplacian matrix of $\overline{G}$, then
$\overline{L}=nI-J-L$, where $I$ is the $n\times n$ identity matrix and $J$ is the $n\times n$ matrix
of all-$1$'s. Consequently, $\uplambda_{\max}\leq n$. 
\item[{\bf 3.}] If $H$ is a (not necessarily induced) subgraph of $G$ on the same vertex-set
$[n]$, and if $\upmu_1\leq\upmu_2\leq\dots\leq\upmu_n$ are the eigenvalues
of the Laplacian of $H$, then $\uplambda_{i}\geq\upmu_{i}$ for all $i\in[n]$. 
\end{itemize}
\end{lem}

Lemma~\ref{lem:tech1} Part {\bf 1}'s proof was discussed during the Introduction (Section~\ref{sec:intro}), and Part {\bf 2} is a straightforward verification, 
but Part {\bf 3} is a more advanced result. 

\section{Largest Eigenvalue of a Comparability Graph.}\label{sec:lecg}

The main goal of this section is to prove the following theorem:
\begin{theo}\label{theo:mainlap}
Let $G=G([n],E)$ be a connected comparability graph with Laplacian matrix $L=L(G)$ and canonical partition $\mathcal{P}=\mathcal{P}(G)$. 
Let $\uplambda_{\max}=\uplambda_{\max}(G)$ be the largest eigenvalue of $L$ and $\eigspace{\max}$ its associated eigenspace. 
Then, the following are true:
\begin{itemize}
\item[{\bf i.}] If $O$ is a transitive orientation of $G$, then: 
$$\dim \left\langle C_O\cap \eigspace{\max}\right\rangle=\dim \left\langle\eigspace{\max}\right\rangle.$$
\item[{\bf ii.}] $\eigspace{\max}\subseteq \bigcup_{O} C_O$, where the union is over all transitive orientations of $G$.
\item[{\bf iii.}] Let $x\in \eigspace{\max}$ be a \emph{u.c.u.v.}. Almost surely: 
\begin{itemize}
\item[{\bf 1.}] If $A\in\mathcal{P}$, then $A$ belongs to a fiber of $x$.
\item[{\bf 2.}] If $A,A'\in\mathcal{P}$ are completely adjacent in $G$, then $A$ and $A'$ belong to
different fibers of $x$.
\item[{\bf 3.}] $x$ induces a transitive orientation of $G^{\mathcal{P}}$. In particular, $G^{\mathcal{P}}$ is a comparability graph. 
\item[{\bf 4.}] All transitive orientations of $G^{\mathcal{P}}$ can be induced by $x$ with positive probability.
\item[{\bf 5.}] If $\upxi$ is a fiber of $x$, then: 
$$G[\upxi]=G[B_1]+\dots +G[B_k],$$ 
where for all $i\in[k]$, $B_i$ is a connected module of $G$ and $G[B_i]$ is a comparability graph. 
\item[{\bf 6.}] $G$ has exactly two transitive orientations if and only if 
$\dim \left\langle\eigspace{\max}\right\rangle=1$ and every fiber of $x$ is an independent set
of $G$.
\end{itemize} 
\item[{\bf iv.}] If $\overline{G}$ is connected, then $\dim \left\langle\eigspace{\max}\right\rangle=1$. If $\overline{G}$ is disconnected, 
then $\dim \left\langle\eigspace{\max}\right\rangle$ is equal to the number of connected components
of $\overline{G}$ minus one. 
\end{itemize}
\end{theo}
\begin{rem}[to Theorem~\ref{theo:mainlap}]\label{rem:mainlap}
In fact, as it will be explained, all transitive orientations of $G$ can be obtained with the following 
procedure: Select an arbitrary transitive
orientation for $G^{\mathcal{P}}$, and select arbitrary transitive orientations for
(the connected components of) each $G[A],A\in\mathcal{P}$ . Therefore, {\bf i-iii} imply an iterative 
algorithm that obtains every
transitive orientation of $G$ with positive probability.  
\end{rem}

The proof of Theorem~\ref{theo:mainlap} will be stepwise and its notation and conventions will carry over to 
the next results, unless otherwise stated. Let us begin with this work.  
 
\begin{prop}\label{prop:cones}
Let $G=G([n],E)$ be a connected comparability graph and let $C_O$ be the (closed convex) cone
corresponding to a transitive orientation $O$ of $G$. Then, $C_O$ contains a non-zero eigenvector of $L$ with eigenvalue
$\uplambda_{\max}$. Furthermore: 
$$\dim \left\langle C_O\cap \eigspace{\max}\right\rangle=\dim \left\langle\eigspace{\max}\right\rangle.$$
\end{prop}
\begin{proof2}
The cases $n=1$ and $n=2$ are easy to verify, so we assume that $n>2$. 

The proof consists of two main steps. Firstly, we will prove that $C_O$ is invariant under left-multiplication
by $L$. Then, we will prove that 
$\dim \left\langle C_O\cap \eigspace{\max}\right\rangle=\dim \left\langle\eigspace{\max}\right\rangle$. 

\begin{flushleft}
\emph{Step 1: $Lx\in C_O$ whenever $x\in C_O$}. 
\end{flushleft}

Take an arbitrary vector $x\in C_O$ and let $\{i,j\}\in E$ with $(i,j)$ in $O$. 
Hence, $x_i\leq x_j$. If we consider the vector $Lx$, then: 
\begin{align*}
(Lx)_j-(Lx)_i&=(x_j\deg j - \sum_{k\in N(j)}x_k)-(x_i\deg i-\sum_{\ell\in N(i)}x_\ell)\\
&=\sum_{k\in N(j)}(x_j-x_k)-\sum_{\ell\in N(i)}(x_i-x_\ell)\\
&=\card{N(i)\cap N(j)}(x_j-x_i)+\sum_{\ell\in N(j)\backslash N(i)}(x_j-x_\ell)\\
&-\sum_{m\in N(i)\backslash N(j)}(x_i-x_m).
\end{align*}
Now, since $O$ is transitive and $G$ is comparability, if $\ell\in N(j)\backslash N(i)$, then 
we must have that $(\ell,j)$ is an edge in $O$, 
so that $x_\ell\leq x_j$ since $x\in C_O$. 
Otherwise, we would require that $\{i,\ell\}\in E$, which is false. Similarly,
if $m\in N(i)\backslash N(j)$, we must have that $(i,m)$ is an edge in $O$, so $x_m\geq x_i$. 
Since also $x_j\geq x_i$ then, we see that
$(Lx)_j-(Lx)_i\geq 0$. Verification of the analogous condition for every edge of $E$ 
shows that indeed $Lx\in C_O$. 

\begin{flushleft}
\emph{Step 2: $\dim \left\langle C_O\cap \eigspace{\max}\right\rangle=\dim \left\langle\eigspace{\max}\right\rangle$.}
\end{flushleft}

Suppose on the contrary
that $\dim \left\langle C_O\cap \eigspace{\max}\right\rangle<\dim \left\langle\eigspace{\max}\right\rangle$. 
Then, there exists $x^{\ast}\in \eigspace{\max}
\big\backslash\text{span}_\R \left\langle C_O\cap \eigspace{\max}\right\rangle$. 
Since $C_O$ is full-dimensional in $\R^{[n]}$, we can write
$x^{\ast}=x-y\text{ for some }x,y\in C_O$, where necessarily either $x\not\in \eigspace{\max}^{\perp}$ or
$y\not\in \eigspace{\max}^{\perp}$. 
In fact, we must have that $x,y\not\in \eigspace{\max}^{\perp}$. Otherwise, 
if $y\in \eigspace{\max}^{\perp}$,
then $x^{\ast}=\displaystyle\lim_{N\rightarrow \infty}L^N(x-y)\big\slash||L^N(x-y)||=\displaystyle\lim_{N\rightarrow \infty}L^Nx\big\slash||L^Nx||\in C_O$
from \emph{Step 1}, and
similarly, if $x\in \eigspace{\max}^{\perp}$
then $x^{\ast}\in -C_O$, so in both cases $x^{\ast}\in \text{span}_\R \left\langle C_O\cap \eigspace{\max}\right\rangle$.
Hence, $0<||L^Nx||,||L^Ny||\leq \uplambda_{\max}^N\max\{||x||,||y||\}$ for
all $N\geq 1$ and, moreover, since both $L^Nx\slash ||L^Nx||$ and $L^Ny\slash||L^Ny||$ 
can be made arbitrarily close to $\text{span}_\R  \left\langle C_O\cap \eigspace{\max}\right\rangle$ 
(in particular, using \emph{Step 1}, each gets close to $C_O\cap \eigspace{\max}$) for large $N$, then 
the same will be true for $\Frac{L^Nx-L^Ny}{\uplambda_{\max}^N\max\{||x||,||y||\}}=c\Frac{L^Nx^{\ast}}{\uplambda_{\max}^N||x^{\ast}||}=cx^{\ast}$, 
where
$c=\frac{||x^{\ast}||}{\max\{||x||,||y||\}}\neq 0$. Therefore, letting $N\rightarrow\infty$, we obtain
that $x^{\ast}\in\text{span}_\R \left\langle C_O\cap \eigspace{\max}\right\rangle$. This contradicts our choice of $x^{\ast}$, so:
$$\eigspace{\max}
\big\backslash\text{span}_\R \left\langle C_O\cap \eigspace{\max}\right\rangle=\emptyset.$$ 
\par
\qed
\end{proof2}

\begin{lem}\label{lem:conmod}
Let $G=G([n],E)$ be a connected comparability graph and let $O$ be a transitive orientation of $G$.  
If $x\in C_O\cap \eigspace{\max}$, $x\neq 0$, satisfies that
$x_u=x_v=\upalpha$ for some $\{u,v\}\in E$ and $\upalpha\in\R$, then there must exist $A\subsetneq [n]$ such that:
\begin{itemize}
\item[{\bf i.}] A is a (proper non-trivial) connected module of $G$ and $u,v\in A$.
\item[{\bf ii.}] $x_i=\upalpha$ for all $i\in A$.
\end{itemize} 
\end{lem}
\begin{proof2}
That such an $x$ may exist is the content of Proposition~\ref{prop:cones}, but we are assuming here that indeed, such an $x$ exists with the stated
properties. 

Consider the maximal (by inclusion) set $A\subseteq [n]$ such that
$G[A]$ is connected, $u,v\in A$, and $x_k=\upalpha$ for all $k\in A$. Primarily, $G[A]$ cannot be equal to $G$, since that would imply
that $x$ is equal to $\upalpha e_{[n]}$, which is impossible. Hence, $G[A]$ is a proper non-trivial connected induced subgraph of $G$. 

We will show that $A$ is a (proper non-trivial connected) module of $G$. 
Suppose on the contrary, that
$A$ is not a module of $G$. Then, there must exist
two vertices $i,j\in A$ such that $N(i)\backslash A\neq N(j)\backslash A$. Consequently, $N(i)\triangle N(j)\backslash A\neq \emptyset$. 
Furthermore, considering a path in $G[A]$ connecting $i$ and $j$, we observe
that we may assume that $i$ and $j$ are adjacent in $G[A]$, so that $\{i,j\}\in E$. Under this assumption, 
suppose now that $(i,j)$ is an edge in $O$. As $O$ is transitive,
we must have that $(i,k)$ is an edge in $O$ whenever $(j,k)$ is. Similarly,
$(k,j)$ must be an edge in $O$ whenever $(k,i)$ is. As such, since $N(i)\backslash A\neq N(j)\backslash A$,
then it must be the case that for $k\in N(i)\triangle N(j)\backslash A$: 
\begin{align*}
\text{If $k\in N(i)$, then $(i,k)$ is an edge in $O$;}\\
\text{and if $k\in N(j)$, then
$(k,j)$ is an edge in $O$.}
\end{align*} 
Left-Multiplying $x$ by the Laplacian of $G$, we obtain:
\begin{gather*}
0=\uplambda_{\max}\upalpha-\uplambda_{\max}\upalpha=\uplambda_{\max}x_j-\uplambda_{\max}x_i\\
=(Lx)_j-(Lx)_i=\displaystyle\sum_{k\in N(j)}(x_j-x_k)-\displaystyle\sum_{\ell\in N(i)}(x_i-x_\ell)\\
=\displaystyle\sum_{k\in N(j)\backslash A\cup N(i)}(x_j-x_k)-\displaystyle\sum_{\ell\in N(i)\backslash A\cup N(j)}(x_i-x_\ell)\\
=\displaystyle\sum_{k\in N(j)\backslash A\cup N(i)}|x_j-x_k|+\displaystyle\sum_{\ell\in N(i)\backslash A\cup N(j)}|x_i-x_\ell|.
\end{gather*}
Since $N(i)\triangle N(j)\backslash A\neq \emptyset$ and $A$ was chosen maximal, 
then at least one of the terms in the last summations must be non-zero and 
we obtain a contradiction. This proves that $A$ is a module of $G$ with the required properties.
\par
\qed
\end{proof2}
\begin{theo}\label{theo:compupo}
Let $G=G([n],E)$ be a connected comparability graph without proper non-trivial connected modules. Then: 
\begin{itemize}
\item[{\bf i.}] Any $x\in \eigspace{\max}\backslash\{0\}$ induces a transitive orientation of $G$.
\item[{\bf ii.}] $\dim \left\langle\eigspace{\max}\right\rangle=1$.
\item[{\bf iii.}] $G$ has exactly two transitive orientations.
\end{itemize} 
\end{theo}
\begin{proof2}
The cases $n=1$ and $n=2$ are easy to check, so we assume that $n>2$.

Fix a transitive orientation $O$ of $G$ and consider the cone $C_O$. Per Proposition~\ref{prop:cones}, we can find at least one 
$x\in C_O\cap \eigspace{\max}$, $x\neq 0$. By Lemma~\ref{lem:conmod} and since $G$ does not have
proper non-trivial connected modules, $x$ must belong to the interior of $C_O$. This establishes {\bf i}. 

To prove {\bf ii}, assume on the contrary, that $\dim \left\langle\eigspace{\max}\right\rangle>1$. 
Consider two \emph{dual} transitive orientations $O$ and $O_{dual}$ of $G$, {\it i.e.} 
$O_{dual}$ is obtained from $O$ by reversion of the orientation
of all the edges. Using {\bf i}, let $y,z\in \eigspace{\max}\backslash\{0\}$ be such that
$y\in \text{int}(C_O)$, $z\in\text{int}(C_{O_{dual}})$, and $z\not\in\text{span}_{\R}\left\langle y\right\rangle$. 
Then, there exists $\upalpha\in(0,1)$
such that $0\neq \upalpha y+(1-\upalpha)z\in \partial \left\langle C_{O}\cap \eigspace{\max}\right\rangle$, contradicting {\bf i}.  

Finally, {\bf iii} follows easily from {\bf i-ii} and Proposition~\ref{prop:cones}.  
\par
\qed
\end{proof2}

The remaining part of the theory will rely heavily on some standard results of 
the spectral theory of the Laplacian (Section~\ref{subsec:stl}). These will be of central importance to establish 
Proposition~\ref{prop:laplamod1}, Proposition~\ref{prop:laplamod2}, and Corollary~\ref{cor:laplamod},
which deal with arbitrary simple graphs. 

\begin{lem}\label{lem:tech2}
Let $G=G([n],E)$ be a complete $p$-partite graph with 
maximal independent sets $A_1,\dots,A_p$. Then, $\uplambda_{\max}=n$ and: 
\begin{align*}
\eigspace{\max}&=\{x\in\R^{\ast[n]}:\text{If $i,j\in A_q$ for some $q\in[p]$, then $x_i=x_j$}\}\\
&=\text{span}_{\R}\left\langle \{e_{A_q}\}_{q\in[p]}\right\rangle\cap\R^{\ast[n]}.
\end{align*}
In particular, $\dim \left\langle\eigspace{\max}\right\rangle=p-1$.
\end{lem}
\begin{proof2}
The complement of $G$ has $p$ connected components, so by  Parts {\bf 1} and {\bf 2} in Lemma~\ref{lem:tech1},
$\uplambda_{\max}=n$ and $\dim  \left\langle\eigspace{\max}\right\rangle=p-1$. Let $b_1,\dots,b_p\in\R$ and let $x\in\R^{\ast[n]}$ be such that 
$x_i=b_q$ for all $i\in A_q$, $q\in[p]$. For any $i\in[n]$, if $i\in A_q$ then
$(Lx)_i=(n-\card{A_q})b_q-(0-\card{A_q}b_q)=nb_q=nx_i$. The set of all such $x$ has dimension $p-1$.   
\par
\qed
\end{proof2}

\begin{lem}\label{lem:tech3}
Let $G=G([n],E)$ be a connected bipartite graph with bipartition
$\{X,Y\}$. Then, $\dim  \left\langle\eigspace{\max}\right\rangle=1$. Furthermore, if $x\in \eigspace{\max}\backslash\{0\}$,
then either $x_i<0$ for all $i\in X$ and $x_j>0$ for all $j\in Y$, or vice-versa.  
\end{lem}
\begin{proof2}
If $G$ is complete $2$-partite, this is a consequence of Lemma \ref{lem:tech2}. Otherwise, as a connected
bipartite graph, $G$ is also a comparability graph and $G$ does not have connected proper non-trivial modules, so Theorem~\ref{theo:compupo} shows that
$\dim  \left\langle\eigspace{\max}\right\rangle=1$ and that $x\in \eigspace{\max}\backslash\{0\}$ induces a transitive orientation of $G$.  
So take $x\in \eigspace{\max}\backslash\{0\}$ and suppose that $x_i=0$, $i\in X$. Then, $(Lx)_i\neq 0$ as
$x$ induces a transitive orientation of $G$ and since $G$ is connected.  
\par
\qed
\end{proof2}
 
We have not found an agreed-upon notation in the literature for the following
objects, so we will need to introduce it here.  

\begin{defn}\label{defn:notation}
Let $G=G([n],E)$ be a simple connected graph, and let $\mathcal{Q}=\{X_1,\dots,X_m\}$ be a partition
of $[n]$ with non-empty blocks.
Then, for all $k\in[m]$:
\begin{itemize}
\item[a.] $G_{X_k}$ will denote the graph on vertex-set $[n]$ and edge-set: 
$$\left\{\{i,j\}\in E:i,j\in X_k\right\}.$$
\item[b.] $R_{X_k}:=\{x\in \R^{\ast[n]}:x_i=0\text{ if }i\not\in X_k,i\in[n]\}$.
\end{itemize}
Also, 
\begin{align*}
R^{\mathcal{Q}}:&=\{x\in\R^{\ast[n]}:x\text{ is constant on each }X_k,k\in[m]\}\\
&=\text{span}_{\R}\left\langle\{e_{X_k}\}_{k\in[m]}\right\rangle\cap \R^{\ast[n]}.
\end{align*} 

\end{defn}

\begin{obs}\label{obs:lap}
In \emph{Definition~\ref{defn:notation}}, 
the linear subspaces $R^{\mathcal{Q}}$ and $R_{X_k}$ for all $k\in[m]$, are mutually orthogonal. 

Furthermore,
any vector $x\in\R^{\ast[n]}$ can be uniquely written as: 
$$x=y+x_1+x_2+\dots+x_m,$$ 
with $y\in R^{\mathcal{Q}}$ and $x_k\in R_{X_k}$, $k\in[m]$. 
\end{obs}

We are now ready to present the results about the space $\eigspace{\max}$ for simple graphs.
Their proofs will use the same language and main ideas, so we will present them contiguously 
to make this resemblance clear. 

\begin{prop}\label{prop:laplamod1}
Let $G=G([n],E)$ be a connected simple graph such that $\overline{G}$ is connected. 
For any fixed proper module $A$ of $G$, the following is true: 
If $x\in \eigspace{\max}$, then $A$ belongs to a fiber of $x$.   
\end{prop}

\begin{prop}\label{prop:laplamod2}
Let $G=G([n],E)$ be a connected simple graph such that $\overline{G}$ is disconnected. Then, $\uplambda_{\max}=n$ and: 
\begin{align*}
\eigspace{\max}=&\{x\in \R^{\ast[n]}: x_i=x_j,\\
&\text{ whenever $i$ and $j$ belong to the same connected component of $\overline{G}$}\}.
\end{align*}
In particular, $\dim \eigspace{\max}$ is equal to the number of connected components of $\overline{G}$ minus one, and 
$G^{\mathcal{P}}$ is a complete $p$-partite graph, where $p$ is the number of connected components of $\overline{G}$.  
\end{prop}

\begin{proof3}[Preliminary Notation for the Proofs of Proposition~\ref{prop:laplamod1} and Proposition~\ref{prop:laplamod2}:]
Let $I$ be the $n\times n$ identity matrix. As usual, $\mathcal{P}=\{A_1,\dots,A_p\}$ will 
be the canonical partition of $G$. 
Let $L$ be the Laplacian matrix of $G$,
$L^{\mathcal{P}}$ be the Laplacian matrix of the copartition subgraph $G^{\mathcal{P}}$ of
$G$, and $L_{A_q}$ be the Laplacian matrix of $G_{A_q}$ for $q\in[p]$.
Firstly, we observe that $L=L^{\mathcal{P}}+\sum_{q=1}^pL_{A_q}$.  
\end{proof3}
\par
\begin{proof2}[Proof of Proposition~\ref{prop:laplamod1}.]
The plan of the proof is to show that the eigenspace of 
$L^{\mathcal{P}}$ corresponding to its largest eigenvalue lives inside
$R^{\mathcal{P}}$, and then to show that this eigenspace is precisely equal to $\eigspace{\max}$.
This will be sufficient since $A\subseteq A_q$ for some $q\in[p]$.

To prove the first claim, first note that left-multiplication by $L^{\mathcal{P}}$ is $R^{\mathcal{P}}$-invariant, 
where the condition that the $A_q$'s are modules is fundamental to prove this. 
Now, for any $x\in\R^{\ast[n]}$, and writing $x=y+x_1+\dots+x_p$ with
$y\in R^{\mathcal{P}}$ and $x_q\in R_{A_q}$, $q\in[p]$, we have that:
$$L^{\mathcal{P}}x=L^{\mathcal{P}}y+\sum_{q=1}^{p}\card{N(A_q)}x_q.$$
Hence,  by Observation~\ref{obs:lap}, if we can show that the largest eigenvalue of $L^{\mathcal{P}}$ is strictly
greater than $\displaystyle\max \{\card{N(A_q)}\}_{q\in[p]}$, the claim will follow. This is what we will do now.

In fact, we will prove that the largest eigenvalue of $L^{\mathcal{P}}$ is strictly greater than 
$\displaystyle\max\{\card{N(A_q)}+|A_q|\}_{q\in[p]}$. 
To check this, first note that both $G^{\mathcal{P}}$ and 
its complement are connected graphs, and that for $q\in[p]$, $A_q$ is both a maximal proper module
and an independent set of $G^{\mathcal{P}}$. For an arbitrary $q\in[p]$, consider the (not necessarily induced) subgraph $H_{\sim q}$
of $G^{\mathcal{P}}$ on vertex-set $A_q\cup N(A_q)$ and whose edge-set is $\left\{\{i,j\}\in E:i\in A_q\text{ and }
j\in N(A_q)\right\}$. Firstly, $H_{\sim q}$ is a complete $2$-partite graph, so its largest eigenvalue is 
precisely $\card{N(A_q)}+|A_q|$ from Lemma~\ref{lem:tech2}. Secondly, since both $G^{\mathcal{P}}$ and its complement are connected, 
there exists a (not necessarily induced)
connected bipartite subgraph $H$ of $G^{\mathcal{P}}$ such that $H_{\sim q}=H[A_q\cup N(A_q)]$ and $H\neq H_{\sim q}$. By Lemma~\ref{lem:tech1} Part {\bf 3} and 
Lemma~\ref{lem:tech3}, the largest eigenvalue of the Laplacian matrix of $H$ must be strictly greater than that of $H_{\sim q}$, 
since any non-zero eigenvector for this
eigenvalue must be non-zero on the vertices of $H$ that are not vertices of $H_{\sim q}$. Also, by the same Lemma~\ref{lem:tech1} Part {\bf 3}, the
largest eigenvalue of $L^{\mathcal{P}}$ must be at least equal to the largest eigenvalue of the Laplacian matrix of $H$. 
This proves the first claim.  

To prove the second claim, note that for $q\in[p]$,
left-multiplication by $L_{A_q}$ is
$R_{A_q}$-invariant. Also, for an arbitrary $x\in\R^{\ast[n]}$ decomposed as above, we have that: 
$$Lx=L^{\mathcal{P}}y+\sum_{q=1}^{p}(\card{N(A_q)}I+L_{A_q})x_q,$$
and this gives the unique decomposition of $Lx$ of
Observation~\ref{obs:lap}. But then, from the proof of the first claim, we note that it suffices to prove that
the largest eigenvalue of $L^{\mathcal{P}}$ is strictly greater than that of
$\card{N(A_q)}I+L_{A_q}$ for any $q\in[p]$. However, from Lemma~\ref{lem:tech1} Part {\bf 1}, we know that
the largest eigenvalue of $L_{A_q}$ is at most $|A_q|$, so the largest eigenvalue
of $\card{N(A_q)}I+L_{A_q}$ is at most $\card{N(A_q)}+|A_q|$. We have already proved that the largest eigenvalue
of $L^{\mathcal{P}}$ is strictly greater than $\max\{\card{N(A_q)}+|A_q|\}_{q\in[p]}$, so the second claim follows. 
\par
\qed
\end{proof2}

\begin{proof2}[Proof of Proposition~\ref{prop:laplamod2}.]
That $G^{\mathcal{P}}$ is a complete $p$-partite graph is clear, so from Lemma~\ref{lem:tech2}, it will suffice to prove that
$\eigspace{\max}$ is exactly equal to the eigenspace of $L^{\mathcal{P}}$ corresponding to its largest eigenvalue 
($=n$). 
This is what we do. 

As in the proof of Proposition~\ref{prop:laplamod1}, we observe that left-multiplication by
$L^{\mathcal{P}}$ is $R^{\mathcal{P}}$-invariant, and that for $q\in[p]$, left-multiplication by 
$L_{A_q}$ is $R_{A_q}$-invariant. For an arbitrary $x\in\R^{\ast[n]}$ with $x=y+x_1+\dots+x_p$,
where $y\in R^{\mathcal{P}}$ and $x_q\in R_{A_q}$, $q\in[p]$, 
and noting that $\card{N(A_q)}=n-|A_q|$ in this case, we have that:
$$Lx=L^{\mathcal{P}}y+\sum_{q=1}^{p}((n-|A_q|)I+L_{A_q})x_q,$$
and this gives the unique decomposition of $Lx$ of
Observation~\ref{obs:lap}. Hence, we will be done if we can show that the largest eigenvalue of 
any of the matrices $L_{A_q}$, $q\in[p]$, is strictly less than $|A_q|$. However, since by construction (from the
definition of canonical partition),
$G[A_q]$ satisfies that its complement is connected, then Lemma~\ref{lem:tech1} Parts {\bf 1} and {\bf 2} imply that
the largest eigenvalue $L_{A_q}$ is strictly less than $|A_q|$, and this holds for all $q\in[p]$. This completes the proof. 
\par
\qed
\end{proof2}
\begin{cor}\label{cor:laplamod}
Let $G=G([n],E)$ be a connected simple graph with canonical partition $\mathcal{P}$ (with $L$ and $\eigspace{\max}$ as usual). 
If $L^{\mathcal{P}}$ 
denotes the Laplacian matrix of $G^{\mathcal{P}}$, then
the eigenspace of $L^{\mathcal{P}}$ corresponding to the largest eigenvalue coincides with $\eigspace{\max}$. 
\end{cor}

Let us now turn back our attention to comparability graphs and to the proofs of Theorem~\ref{theo:lowlap} and
Theorem~\ref{theo:mainlap}. Comparability graphs are, as anticipated, specially amenable to
apply the previous two propositions and their corollary. In fact, the following result 
already establishes most of Theorem~\ref{theo:mainlap}.
\begin{prop}\label{prop:laplamodcomp}
Let $G=G([n],E)$ be a connected comparability graph with canonical partition
$\mathcal{P}$. 
\begin{itemize}
\item[{\bf i.}] For $x\in\eigspace{\max}$ a \emph{u.c.u.v.}, the following hold true almost surely:
\begin{itemize}
\item[{\bf 1.}] If $A\in\mathcal{P}$, then $A$ belongs to a fiber of $x$.
\item[{\bf 2.}] If $A,A'\in\mathcal{P}$ are completely adjacent in $G$, then $A$ and $A'$ belong to
different fibers of $x$.
\item[{\bf 3.}] $x$ induces a transitive orientation of $G^{\mathcal{P}}$. 
In particular, $G^{\mathcal{P}}$ is a comparability graph. 
\item[{\bf 4.}] If $\upxi$ is a fiber of $x$, then: 
$$G[\upxi]=G[B_1]+\dots +G[B_k],$$ 
where for all $i\in k$, $B_i$ is a connected module of $G$ and $G[B_i]$ is a comparability graph. 
\end{itemize} 
\item[{\bf ii.}] If $\overline{G}$ is connected, then $\dim \left\langle\eigspace{\max}\right\rangle=1$. Also, $G^{\mathcal{P}}$ has exactly two transitive orientations and each 
can be obtained with probability $\frac{1}{2}$ in {\bf i}.
\item[{\bf iii.}] If $\overline{G}$ is disconnected, then $\dim \left\langle\eigspace{\max}\right\rangle=p-1$, where 
$p$ is the number of connected components of $\overline{G}$. Also, $G^{\mathcal{P}}$ has exactly $p!$ transitive orientations and each can be
obtained with positive probability in {\bf i}. 
\end{itemize} 
\end{prop}
\begin{proof2}
We will work on each case, whether $\overline{G}$ is connected or disconnected, separately.

\emph{Case 1: $\overline{G}$ is connected.}

From Proposition~\ref{prop:cones}, take any $x\in C_O\cap \eigspace{\max}$,  $x\neq 0$, for some transitive orientation $O$
of $G$. From Proposition~\ref{prop:laplamod1}, we know that $x$ is constant on each $A\in\mathcal{P}$, so {\bf i.1} holds. 
Moreover, since  
the elements of $\mathcal{P}$ are the maximal proper modules of $G$, then Lemma~\ref{lem:conmod} shows that for 
completely adjacent $A,A'\in\mathcal{P}$, $x_i\neq x_j$ whenever $i\in A$ and $j\in A'$, so {\bf i.2} holds. 
Now, since the orientation of $G^{\mathcal{P}}$ induced by $x$ is then equal to the
restriction of $O$ to the edges of $G^{\mathcal{P}}$, we observe that for $A,A'$ as above,  
the edges $\left\{\{i,j\}\in E:i\in A\text{ and }j\in A'\right\}$ are oriented in $O$ in the \emph{same direction} (either from $A$ to $A'$, or vice-versa). 
Since $O$ is transitive, this immediately
implies that its restriction to $G^{\mathcal{P}}$ is transitive, so $G^{\mathcal{P}}$ is a comparability graph and 
{\bf i.3} holds. 
Notably, this holds for any choice of $O$. If $\upxi$ is a fiber of $x$, then we can write $G[\upxi]$ as a disjoint union of its connected components,
say $G[\upxi]=G[B_1]+\dots+G[B_k]$. On the one hand, the restriction of $O$ to any induced subgraph of $G$ is transitive, so
$G[\upxi]$ is a comparability graph, and also each of its connected components. On the other hand, from {\bf i.2}, each $B_i$ with $i\in[k]$ satisfies that
$B_i\subseteq A$ for some $A\in\mathcal{P}$, and moreover, $G[B_i]$ is a connected component of $G[A]$, so $B_i$ is a module
$G$ since $B_i$ is a module of $A$ and $A$ is a module of $G$. This proves {\bf i.4}. 

As $G^{\mathcal{P}}$ does not have proper non-trivial connected modules, 
from Theorem~\ref{theo:compupo} and
Corollary~\ref{cor:laplamod}, we obtain that $\dim \left\langle\eigspace{\max}\right\rangle=1$. Also, $G^{\mathcal{P}}$ has exactly two transitive orientations and
each can be obtained with probability $\frac{1}{2}$ from $x\in\eigspace{\max}$ a \emph{u.c.u.v.}, proving
{\bf ii}. 

\emph{Note:} In fact, then, it follows that for any $x\in \eigspace{\max}\backslash\{0\}$, necessarily $x\in C_O$ or $x\in\mathcal{C}_{O_{dual}}$, 
where $O$ is the orientation used in the proof, and 
$O_{dual}$ is the dual orientation to $O$.  

\emph{Case 2: $\overline{G}$ is disconnected.}

This is precisely the setting of Proposition~\ref{prop:laplamod2}, so {\bf i.1-3} and {\bf iii} follow after
noting that, firstly, $p$-partite graphs are comparability graphs, and secondly, their
transitive orientations are exactly the acyclic orientations of their edges such that: 
\begin{itemize}
\item[] For every pair of maximal independent sets, all the edges between them (or having endpoints on both sets), are oriented
in the same direction. 
\end{itemize}
The proof
of ${\bf i.4}$ goes exactly as in \emph{Case 1}. 
\par
\qed
\end{proof2}

\begin{cor}\label{cor:transcomp}
Let $G=G([n],E)$ be a connected comparability graph with
canonical partition $\mathcal{P}$, and let $O$ be a transitive orientation of $G$. 
Then, {\bf (1) } the restriction of $O$ to each of $G^{\mathcal{P}}$ and 
$G[A],A\in\mathcal{P}$, is transitive. 

Conversely, {\bf (2)} if we select arbitrary transitive orientations for each of 
$G^{\mathcal{P}}$ and $G[A],A\in\mathcal{P}$, and then take the union of these, we obtain a transitive orientation for
$G$. 
\end{cor}
\begin{proof2}
Statement {\bf (1)} follows from Proposition~\ref{prop:laplamodcomp} and Proposition~\ref{prop:cones},
since $\dim \left\langle C_O\cap \eigspace{\max}\right\rangle=\dim \left\langle\eigspace{\max}\right\rangle$. 

For {\bf (2)}, select transitive orientations for each of $G^{\mathcal{P}}$ and $G[A],A\in\mathcal{P}$, and let
$O$ be the orientation of $E$ so obtained. Since each element of $\mathcal{P}$ is independent in
$G^{\mathcal{P}}$ and since the restriction of $O$ to $G^{\mathcal{P}}$ is transitive, then: 
\begin{itemize}
\item[($\star$)] For $A,A'\in\mathcal{P}$ completely adjacent, the edges between
$A$ and $A'$ must be oriented in $O$ in the same direction. 
\end{itemize}
This rules out
the existence of directed cycles in $O$, so $O$ is acyclic. Now, if $O$ is not transitive,
then there must exist $i,j,k\in[n]$ such that $(i,j)$ and $(j,k)$ are in $O$ but not
$(i,k)$. By the choice of $O$, it must be the case that exactly two among $i,j,k$ 
belong to the same $A\in\mathcal{P}$, and the other one to a different $A'\in\mathcal{P}$. 
The former cannot be $i$ and $k$, per the argument above ($\star$). Hence, without loss of generality,
we can assume that $i,j\in A$ and $k\in A'$. But then, $A$ and $A'$ must be completely adjacent and
$(i,k)$ must exist in $O$, so we obtain a contradiction. 

\emph{Note:} The argument for {\bf (2)} is essentially found in \cite{pg}.    
\par
\qed
\end{proof2}

\begin{cor}\label{cor:compnonupo}
Let $G=G([n],E)$ be a connected comparability graph with at least one proper non-trivial connected module $B$,
and canonical partition $\mathcal{P}$. Then, 
$G$ has more than two transitive orientations.
\end{cor}
\begin{proof2}
Suppose, on the contrary, that $G$ has only two transitive orientations. 
We will prove that, then, $G$ cannot have proper non-trivial connected modules and
so $B$ does not exist. 

From Corollary~\ref{cor:transcomp}
and Proposition~\ref{prop:laplamodcomp}.{\bf ii-iii}, a necessary condition for $G$
to have no more than two transitive orientations is:
\begin{itemize}\label{itm:cond}
\item[($\star$)] $G=G^{\mathcal{P}}$, and either $\overline{G}$ is connected or it has exactly two connected components.
\end{itemize}
Now, if $\overline{G}$ is connected, then $B\subseteq A$ for some
$A\in\mathcal{P}$ by
Corollary~\ref{cor:modbasic}, so $B$ is an independent set of $G$ since
$A$ is independent. This contradicts the choice of $B$. Also, if $\overline{G}$ has two connected components, then
$G$ is a complete bipartite graph. However, it is clear that no such $B$ can exist in a 
complete bipartite graph.
\par
\qed
\end{proof2}

\begin{proof2}[Proof of Theorem~\ref{theo:mainlap}.]
The different numerals of this result have, for the most part, already been proved.
\begin{itemize}
\item[-] {\bf i} was proved in Proposition~\ref{prop:cones}.
\item[-] {\bf ii} was proved in Proposition~\ref{prop:laplamodcomp} for
the case when $\overline{G}$ is connected (See \emph{Note}). In the general case, {\bf ii} follows from
Proposition~\ref{prop:laplamodcomp}.{\bf i.1-3} and Corollary~\ref{cor:transcomp} Statement {\bf (2)} for $x\in\eigspace{\max}$ a \emph{u.c.u.v.}, and then for all
$x\in\eigspace{\max}$ since the cones $C_O$ (with $O$ an acyclic orientation of $E$) are closed.
\item[-] {\bf iii.1-5} and {\bf iv} are precisely Proposition~\ref{prop:laplamodcomp}.
\item[-] For {\bf iii.6}, from Corollary~\ref{cor:compnonupo} and Theorem~\ref{theo:compupo}.{\bf iii}, $G$ has exactly two transitive orientations
if and only if $G$ has no proper non-trivial connected modules. Now,
if $G$ has no proper non-trivial connected modules, then Proposition~\ref{prop:laplamodcomp}.{\bf i.4} shows that 
the fibers of $x$ are independent sets of $G$ and Theorem~\ref{theo:compupo}.{\bf ii} gives
$\dim \left\langle\eigspace{\max}\right\rangle=1$. Conversely, if the fibers of $x$ are independent sets of $G$,
then $G=G^{\mathcal{P}}$. Furthermore, per Proposition~\ref{prop:laplamodcomp}.{\bf ii-iii},  
if $\dim \left\langle\eigspace{\max}\right\rangle=1$, then $\overline{G}$ has at most two connected components. 
Hence, $G=G^{\mathcal{P}}$ and $\overline{G}$ has at most two connected components,
so we obtain precisely the setting of ($\star$) in Corollary~\ref{cor:compnonupo}. 
Consequently, $G$ cannot have proper non-trivial connected modules.  
\end{itemize}
\par
\qed
\end{proof2}

\section{A characterization of comparability graphs.}\label{sec:ccg} 

This section offers a curious novel characterization of comparability graphs that results from our 
theory in Section~\ref{sec:lecg}.

\begin{theo}\label{theo:lowlap}
Let $G=G([n],E)$ be a simple undirected graph with Laplacian matrix $L$, and let $I$ be the
$n\times n$ identity matrix.

Then, $G$ is a comparability graph if and only if there exists 
$\upalpha\in\N$ and an acyclic orientation $O$ of $E$, such that $C_O$ is invariant under left-multiplication by
$\upalpha I+L$.

If $G$ is a comparability graph, the orientations that satisfy the condition are precisely the transitive 
orientations of $G$, and we can take $\upalpha=0$ for them.  
\end{theo}
\begin{proof2}
If $G$ is a comparability graph and $O$ is a transitive orientation of $G$, then
\emph{Step 1} of Proposition~\ref{prop:cones} shows that indeed, $Lx\in C_O$
whenever $x\in C_O$. Clearly then, for all $\upalpha\in\R_{\geq 0}$, $(\upalpha I+L)x\in C_O$
whenever $x\in C_O$. 

Suppose now that $G$ is an arbitrary simple graph, and let 
$O$ be an acyclic orientation (of $E$) that is not a transitive orientation of
$G$. Then, there exist $i,j,k\in[n]$ such that $(i,j)$ and $(j,k)$ are in $O$ but not
$(i,k)$, and the following set is non-empty: 
\begin{gather*}
X:=\{k\in[n]:\text{ there exist }i,j\in[n]\text{ and directed edges }\\
(i,j),(j,k)\text{ in }O\text{, but }(i,k)\text{ is not in }O\}.
\end{gather*}
In the partial order on $[n]$ induced by $O$, take some $\ell\in X$ maximal, and
consider the principal order filter $\ell^{\vee}$ whose unique minimal element is $\ell$.
The indicator vector of $\ell^{\vee}$ is $e_{\ell^{\vee}}$.
Then, $e_{\ell^{\vee}}\in C_O$. Now, choose $i,j\in[n]$ so that 
$(i,j)$ and $(j,\ell)$ are in $O$ but not $(i,\ell)$. As $\ell$ was chosen maximal in $X$, for every
$k\in \ell^{\vee}, k\neq \ell$, then both $(i,k)$ and $(j,k)$ are in $O$. Therefore, we have:
\begin{align*}
(Le_{\ell^{\vee}})_i&=-\card{\ell^{\vee}}+1,\text{ and }\\
(Le_{\ell^{\vee}})_j&=-\card{\ell^{\vee}}.
\end{align*}
Hence, $(Le_{\ell^{\vee}})_i>(Le_{\ell^{\vee}})_j$ and $Le_{\ell^{\vee}}\not\in C_O$ since
$(i,j)$ is in $O$. Since actually $e_{\ell^{\vee}}\in \partial C_O$, then
$(\upalpha I+L)e_{\ell^{\vee}}\not\in C_O$ for $\upalpha\in\R_{\geq 0}$.
\par
\qed
\end{proof2}

\nocite{*}
\bibliographystyle{abbrvnat}
\bibliography{tt}
\label{sec:biblio}

\end{document}